\documentclass[11pt]{article}
\usepackage{geometry}                
\usepackage{graphicx}
\usepackage{amssymb}
\usepackage{epstopdf}
\usepackage{amsmath}
\usepackage{comment}
\usepackage{xcolor}
\usepackage[title]{appendix}
\usepackage{ulem}

\newcommand{\real}{{\mathbb R}}
\newcommand{\OO}{\mathcal{O}}

\newcommand{\nn}{\mathbf{n}}
\newcommand{\tnn}{\tilde{\mathbf{n}}}
\newcommand{\vv}{\mathbf{v}}
\newcommand{\tvv}{\tilde{\mathbf{v}}}
\newcommand{\talpha}{\tilde{\alpha}}

\newcommand{\Lu}{L_{u^*}}

\newtheorem{theorem}{Theorem}

\newtheorem{remark}{Remark}

\begin{document}

\title{Damped and Driven Breathers and Metastability}

\author{Daniel A. Caballero \and C. Eugene Wayne}


\maketitle

\begin{abstract}  
{ In this article we prove the existence of a new family of periodic solutions for discrete,
nonlinear Schr\"odinger equations subject to spatially localized driving and damping.  They provide an alternate description of the metastable
behavior in such lattice systems which agrees with
previous predictions for the evolution of metastable states  while providing more
accurate approximations to these states. We analyze the stability of these breathers, finding a very small positive eigenvalue whose eigenvector lies almost tangent to the surface of the cylinder formed by the family of breathers. This causes solutions to slide along the cylinder without leaving its neighborhood for very long times.
}
\end{abstract}

\section{Introduction and background}

Metastable systems are common throughout multiple disciplines including biology, physics, machine learning and chemistry.  In this paper we reexamine the existence and properties of metastable states in (finite) lattices of coupled nonlinear oscillators.  We are specifically interested in a finite discrete nonlinear Schr\"odinger equation (DNLSE) chain in which energy is forced to flow through the system from one end to the other.  Stochastically perturbed systems of this type have been used as models to
understand the microscopic origins of heat transport and steady state behavior in
non-equilibrium statistical systems \cite{Eckmann:1999}.  A key question
in such studies is how long it takes for the system to converge to its invariant measure.  This controls
how long it takes for the system to exhibit ergodic properties - in some cases making that time scale so long that these properties are effectively unobservable \cite{Iubini:2020} \cite{Danieli:2017} \cite{Danieli:2019} \cite{Hairer:2009} {\cite{Pace2019}.  While statistical mechanics
methods have been applied to such systems, it has been been observed that the presence of breathers in the system prevents the definition
of a well-defined thermodynamic temperature for the state.
\cite{Rasmussen1999} \cite{Levy2018}.}   In both numerical studies and theoretical investigations, it has been observed that the solution can be ``trapped'' for long times in small regions of the phase space.  This long-time trapping
of the solution in restricted regions of phase space is the metastability we will be studying, and in
this paper we look for geometric structures  in the phase space, and an analytical mechanism that can help explain this trapping.

In recent work, it has been proven that for chains of oscillators subject to very weak, localized damping, there are open sets of initial data in which this metastable behavior
can be approximated for very long times by ``breathers'', i.e. temporally periodic, spatially localized solutions of the {\em undamped} equation \cite{Eckmann:2018} \cite{Flach:2008} \cite{Eckmann:2020}.  The proof used a modulation approach in which the solution of the damped system was decomposed into a breather whose frequency and phase were allowed to vary (slowly) in time, and a correction term which was proven to remain small for very long times. A similar effect is present in the experimental technique of self-localization, where damping is turned on at both ends of a dNLSE chain for a short period of time to prepare the system into a localized state after which the dissipation is turned off \cite{Hennig2013} \cite{Livi2006}. Thus, we expect that while employing this technique, the system is in fact in a metastable state similar to the one studied above and it is this metastability which causes the phenomenon of self-localization.

In the present paper we prove the existence of a new family of breather type solutions 
which give even better approximations to the metastable states. { We are specifically
interested in modeling the types of states that have been found to inhibit the
convergence of models of heat conduction toward their invariant measure
 \cite{Danieli:2019} \cite{Hairer:2009} \cite{Cuneo:2015}.  These states are ones
 in which the motion of the system is approximately periodic, and there is a strong localization
 effect of the system. In additional the dissipation in the system typically only occurs near the ends,
 either through an explicit dissipative term, or through the coupling of the end of the 
 system to a heat bath.  This localization of the dissipation is important, and as pointed
 out in  \cite{Hairer:2009}(p.1001) is quite different from models in which the dissipation
 acts equally on all sites of the system.}  {This new family of breather type solution is} constructed,
not by regarding the damped system as a perturbation of the undamped case, but instead by
adding an additional {\em very small} perturbation to the damped system.  We prove that this perturbed system has breather solutions, {and we then analyze the stability of these solutions.
We show that these new localized
states (which we dub {\em damped and driven breathers}) reproduce previous predictions for
the metastable behavior, but also provide more accurate approximations to the metastable
states themselves.  An additional advantage of these solutions is that they allow us to use standard techniques for analyzing the stability of fixed points in studying the metastable states.
}

{
There has been an immense amount of work on breathers in general nonlinear lattices,
including the addition of damping and driving terms.  Such systems arise naturally
in optical cavities, for example, \cite{Peschel:2004}, \cite{Prilepsky:2012}.  The existence
and properties of the breathers in such systems have been studied both numerically and analytically
in these papers and others such as \cite{Efremidis:2003}, \cite{Hennig:1999}, \cite{Marin:2001}.  However, 
in these models the dissipation is spread throughout the system 
and this  creates
very different issues and effects, both physically
and mathematically.  Indeed, one of the first rigorous investigations of breathers
included the possibility of dissipation \cite{MacKay:1997}, but again, the types of systems
considered in that work do not include those with localized dissipation of the type studied here.

Much closer to our work is that of Khomeriki \cite{Khomeriki:2004}, and Maniadis, et al \cite{Maniadis:2006}.
In the first of these papers, the author studies a discrete NLS equation with driving and one end
and damping at the other, very similar to the model we study.  The author finds that 
for sufficiently large driving, one can create moving solitons in band gaps for this model.  For
small driving amplitudes one finds stationary breathers similar to those analyzed here, but the
focus of  \cite{Khomeriki:2004} is not on the sort of metastability issues considered in this
paper nor is the stability of those solutions analyzed.  Maniadis, et al, also consider a chain of
nonlinear Hamiltonian oscillators driven at one end, though in this case they focus on a
discrete Klein-Gordon chain.   They also find numerically that one can generate a large
family of moving breather type solutions in this model by varying the driving force,
and they  track the types of solutions generated by a numerically implemented bifurcation analysis.
Again, the focus of this paper is not of the type of metastable behavior that we are interested in.

Finally, one very interesting possible connection of our work to current research that was pointed
out by a referee is to PT-symmetry breaking in nonlinear lattices \cite{Konotop:2016}.
If one chooses the strength of the damping and driving to be equal, then the linearization
of our problem satisfies the PT symmetry property.  Quantum mechanical systems with this
type of perturbation present at two sites as a sort of ``impurity''  have been studied
in \cite{Jin:2009} and \cite{Joglekar:2010}.  There, the authors find that if the strength
of the perturbation is sufficiently large, one has PT-symmetry breaking and the appearance
of exponentially growing or decaying states.  One could then ask whether there is a bifurcation
that occurs when the driving and damping differ from one another and whether or not this
could be used to give an alternate proof of the existence of the sort of states we study.
Since in our particular situation, the values of the driving and damping differ by several
orders of magnitude, we haven't pursued this approach here, but it seems like an interesting
question for future investigation.

Impurity states have also been studied in systems without PT-symmetry (see for example,
\cite{Theocharis:2009}) where they have again been shown to give rise to localized breathing
modes.  However, the mechanism in this case is rather different than the one we observe since
it arises from exponentially localized eigenstates of the linearized problem created by the
potential well associated with the impurity, whereas the states we study are perturbations 
of inherently nonlinear states.

We conclude the introduction with an overview of the remainder of the paper.  Our goal is
to understand the trapping of energy in localized regions of the lattice which is 
observed in studies of heat transport in coupled lattice of nonlinear Hamiltonian oscillators.
These ``trapped'' regions correspond to states in which the energy is highly localized in
particular spatial regions, and exhibits a very slow dissipation rate.  With this goal,
we first prove the existence of a class of spatially localized, temporally periodic
solutions in a system subject to weak damping and driving.  We then derive
precise estimates on these states which we show reproduce aspects of the 
metastable or trapped states observed in the numerical studies of such systems.
The states we construct occur in continuous families and we then analyze the stability of this
family.  { In particular, we prove that the linearization of the system at each of these
fixed points has one zero eigenvalue, corresponding to a symmetry of the system, one very small positive eigenvalue, and
all the remaining eigenvalues have strictly negative real part.   We then argue that these stability properties explain why solutions are trapped near
these states for long periods - in particular, the eigenvalues with negative real part imply that nearby solutions are attracted
toward the family of fixed points, while the small positive eigenvalue leads to a very slow drift along
that family.}  Finally, we present numerics which show that 
solutions which start close to some member of one of these families of states
do remain close to that family for very long times.

}

\subsection{Past results}

Following \cite{Eckmann:2018}\, \cite{Eckmann:2020}, we work with a very
specific system of coupled oscillators, namely the discrete nonlinear
Schr\"odinger (dNLS) equation:

\begin{equation}\label{eq:NLS}
-i \dot{z}_j = -\epsilon (\Delta z)_j + |z_j|^2 z_j \ ,\ \ j=1, 2, \dots , N\ .
\end{equation}
Note that if $\epsilon=0$, this system becomes $N$, uncoupled, nonlinear
oscillators, and we then have trivial, localized, periodic solutions
in which one site rotates with non-zero angular frequency, and all
other sites have zero amplitude.  This is sometimes referred to as the 
{\it anti-integrable limit.}

We add a driving term to the first site of our system, and a damping term to the last site by
modifying the equation as
\begin{equation}\label{eq:DDNLS}
-i \dot{z}_j = -\epsilon (\Delta z)_j + |z_j|^2 z_j - i \beta \delta_{j,1} z_1 + i \gamma \delta_{j,N} z_N \ ,\ \ j=1, 2, \dots , N\ .
\end{equation}
The addition of damping and driving to the equation allows one to
study energy flow from one end of the system to the other.
The existence of breathers in the dNLS equation with localized forcing but distributed damping was studied in \cite{Panayotaros:2014}. Throughout the paper we will be referring to the $\gamma=0,\beta=0$ case as undamped, $\gamma\ne0,\beta=0$ as damped, and $\gamma\ne0,\beta\ne0$ as damped and driven.
{  We emphasize again that the choice of our particular type of damping and driving 
is motivated by the desire to understand the very slow evolution of models of heat conduction
toward their equilibrium state and on the basis of numerical experiments the most slowly
decaying states appear to be those in which the energy is localized as far as possible from
the sites at which the damping occurs.  However, we believe that if one chose to impose the
damping and driving at other sites, one would also observe localized periodic solutions of the type
we construct - however they would correspond to more rapidly decaying states in the
heat conduction models.}

If we now translate to a rotating frame of reference $z_j(t) = e^{i \omega t}\zeta_j$, then
\begin{equation}\label{eq:NLS_rotating}
-i \dot{\zeta}_j = -\epsilon (\Delta \zeta)_j  - \omega \zeta_j + | \zeta_j |^2 \zeta_j\
- i \beta \delta_{j,1} \zeta_1 + i \gamma \delta_{j,N} \zeta_N  .
\end{equation}

Finally, we will often want to work with real variables so we decompose $\zeta_j = p_j + i q_j$. 
Then \eqref{eq:NLS_rotating} can be written as the system of differential equations:
\begin{eqnarray}\label{eq:NLS_P}
\dot{p}_j &=& \epsilon (\Delta q)_j +\omega q_j - (|p_j |^2 + |q_j |^2)q_j \\ \nonumber
&& \qquad \qquad \qquad \qquad \qquad + \beta \delta_{j,1} p_1
- \delta_{j,N} p_N  \\ \label{eq:NLS_Q}
\dot{q}_j &=& -\epsilon (\Delta p)_j -\omega p_j +(|p_j |^2 + |q_j |^2)p_j \\ \nonumber
&& \qquad \qquad \qquad \qquad \qquad + \beta \delta_{j,1} q_1
- \delta_{j,N} q_N  \ .
\end{eqnarray}
Past research on systems of this type has focused mainly on the Hamiltonian case in which the damping
and driving are both absent - i.e. $\gamma = \beta =0$.  In this case, there is an extensive
theory investigating the existence and properties of ``breathers''.  These spatially localized
but temporally periodic solutions  are known to exist for very general types of lattices of
nonlinear Hamiltonian oscillators including the NLS equation, \cite{MacKay:1994} \cite{Flach:1998} \cite{Kevrekidis:2001}.

For the particular case of \eqref{eq:NLS_P}-\eqref{eq:NLS_Q}, one has the following result:

\begin{theorem}\label{th:undamped_breathers} (see \cite{Eckmann:2020}; Theorem 1.6)
If $\beta = \gamma=0$, there exists $\epsilon_0>0$ and $\Delta \Omega > 0$ such that
for all $|\epsilon | < \epsilon_0$, and $|\omega -1 | < \Delta \Omega$, the equations \eqref{eq:NLS_P}-\eqref{eq:NLS_Q}
have a family of fixed points $(p_j^*(\epsilon,\omega),0)$ which vary smoothly with $\epsilon$ and $\omega$,
and such that $p_1^* = 1+{\cal{O}}(\epsilon,| 1  -\omega|)$, $p_j^* = {\cal{O}}(\epsilon^{j-1})$, $j=2, \dots , n$,
and $q_j^* = 0$.
\end{theorem}

\begin{remark} Note that when we say that \eqref{eq:NLS_P}-\eqref{eq:NLS_Q}
has a family of fixed points of the form $(p_j^*(\epsilon,\omega),0)$,
we mean that the $q_j$ components of the fixed point are zero.
\end{remark}

\begin{remark} Due to the invariance of \eqref{eq:NLS_rotating} under complex rotations
$\zeta_j \to e^{i \theta} \zeta_j$, each of the fixed points constructed in Theorem \ref{th:undamped_breathers}
corresponds to a circle of fixed points. 
\end{remark}

\begin{remark} Because the fixed points in Theorem  \ref{th:undamped_breathers} are constructed with
the aid of the Implicit Function Theorem, for each $(\epsilon,\omega)$, they are the unique fixed points in
a neighborhood of $p_1=1, p_2=p_3=\dots =p_N=0$ (up to the complex rotation noted in the previous
remark).
\end{remark}

\begin{remark} Note that each of the fixed points constructed in Theorem  \ref{th:undamped_breathers}
corresponds to a periodic solution with angular frequency $\omega$ for \eqref{eq:DDNLS}.  If we fix
$\epsilon$ and consider the family of solutions as a function of $\omega$, combined with the circle of solutions given by the complex rotations we see that in the phase space of \eqref{eq:DDNLS} we have a cylinder
filled with periodic orbits.
\end{remark}

\section{Existence and properties of damped and driven breather for NLS}

\subsection{The implicit function theorem and the existence of breathers for the damped and driven system}

{

Coupled systems of nonlinear oscillators like \eqref{eq:DDNLS} have often been studied as models for heat transport in solid state matter.  In these studies, the first and last oscillator are coupled to stochastic heat baths with different  temperatures, to force an energy flow through the system, and a damping term is added to insure that the total energy in the system does not grow without bound.  It has been observed that the convergence to the invariant
measure in such systems is often extremely slow due to
the presence of metastable states in the phase space.
In prior work,  the energy flow through the system is enforced by placing a localized damping term at one end of the system and placing all the initial energy of the system at the opposite end
\cite{Cuneo:2017} \cite{Eckmann:2018}.  One can prove that the energy in the system tends to zero as $t\to \infty$ and hence it must flow through the system in order to dissipate in this way.  In \cite{Eckmann:2020} it is proven that for very long times during this dissipative process the system remains in a neighborhood of the cylinder of breather solutions for the undamped system.

In this paper we argue that a new family of solutions for a system subject to both damping and (very weak) driving  provides even better approximations to the metastable states than do the undamped breathers.  The existence of this family of breathers follows from the following theorem.

\begin{theorem}\label{th:last_phase}  For $\epsilon$  and $\gamma$ sufficiently small
and all $\omega$ in a neighborhood of $\omega =1$, 
there exists a unique value of $\beta$ such that the equations \eqref{eq:NLS_P}-\eqref{eq:NLS_Q} have a stationary
solution in a neighborhood of the undamped breather.  This family of fixed points depends smoothly
on $\gamma$, $\omega$ and $\epsilon$.
\end{theorem}

\begin{remark} Note that stationary solutions of \eqref{eq:NLS_P}-\eqref{eq:NLS_Q} correspond to periodic solutions with angular frequency $\omega$ of \eqref{eq:DDNLS}.
\end{remark}

The proof of this theorem follows from an application of the Implicit Function Thoerem.  Details are provided in Appendix \ref{ap:IFT}.

}
\subsection{Approximating the damped and driven breathers}\label{sec:approx}

{

We wish to show that the damped and driven breathers constructed in the previous subsection provide a better approximation to the metastable states that occur during energy dissipation in these systems than do the undamped breathers.  To do so, we need accurate approximations to these solutions which we derive in this subsection.  Because we proved the existence of these breathers with the aid of
implicit function theorem, we know  that they are analytic functions of the parameters 
$\epsilon$, $\gamma$ and $\omega$, and hence we can derive power series approximations
for them.  We describe those approximations in this section and provide the details of
their derivation in Appendix \ref{ap:IFT}.

}
In Appendix \ref{ap:IFT}, we define functions $f_j(p_2, p_3, \dots p_N, q_2, q_3, \dots, q_N; \epsilon, \omega, \gamma, p_1,q_1)$ whose
solutions give the breather solutions.  For example, 
$$ f_2(p_2,p_3,\dots,p_N,q_2,q_3,\dots,q_N;\epsilon,\omega,\gamma,p_1,q_1) = -\epsilon(p_{3}+p_{1}-p_2) - \omega p_2 + (p_2^2+q_2^2)p_2$$
 - i.e.
the component of the vector field corresponding to $\dot{q}_2$.

In order to derive approximations to the breathers, we consider
\begin{eqnarray}
 f_j(p_2,p_3,\dots,p_N,q_2,q_3,\dots,q_N;\epsilon,\omega,\gamma,p_1,q_1) &=& 0 \\
  f_{j+N}(p_2,p_3,\dots,p_N,q_2,q_3,\dots,q_N;\epsilon,\omega,\gamma,p_1,q_1) &=& 0\ ,
\end{eqnarray}
for $j=2, \dots , N$.  We want to approximate the solutions of these equations when $\epsilon$
and $\gamma$ are small and $p_1$ and $\omega$ are near $1$.

Note the equations for $f_2,  \dots f_{N-1}, f_{N+2}, \dots ,f_{2N-1}$  do not involve $\gamma$ or $\beta$,
and hence are the same as those for the undamped
breathers.  Thus we can construct power-series approximations for these
equations exactly as in the undamped case - i.e., 
$p_2^*= -\frac{\epsilon p_1}{\omega} + \dots$, $p_3^*= \frac{\epsilon p_2^*}{\omega} 
= (\frac{\epsilon}{\omega})^2 p_1+ \dots $, $\dots$, $p_{N-1}^* =  (-\frac{\epsilon}{\omega})^{N-2} p_1
+ \dots$, while all of the $q_j^*$ are zero to this order in $\epsilon$.
(See Section 2 of \cite{Eckmann:2018} for more details.) Note that because the
implicit function theorem guarantees that the solution is analytic and unique, if we can find a 
consistent power series approximation to the solution to a given order in $\epsilon$, $\gamma$, $\dots$,
then it must correspond to the actual solution.  In order to find the leading order terms
in $p_N^*$ and $q_N^*$, we must look simultaneously at the equations for $f_N$ and $f_{2N}$.
We find that (to lowest order in $\epsilon$ and $\gamma$) they satisfy:
\begin{eqnarray}
\omega q_N^* - \gamma p_N^* &=& 0 \\ \nonumber
\epsilon p_{N-1}^* - \omega p_N^* - \gamma q_N^* &=& 0\ .
\end{eqnarray}
Inserting our leading order expansion for $p_{N-1}^*$, and inverting the matrix $\left( \begin{array}{cc}
\omega & -\gamma \\ - \gamma & - \omega \end{array}\right)$, we find:
\begin{eqnarray}\label{eq:q*p*} \nonumber
\left( \begin{array}{c} q_N^* \\ p_N^* \end{array} \right) &=& \frac{1}{\omega^2 + \gamma^2} 
\left( \begin{array}{cc} \omega & -\gamma \\ -\gamma &  -\omega \end{array} \right)
\left( \begin{array}{c} 0 \\ \frac{(-1)^{N-2} \epsilon^{N-1}}{\omega^{N-2}} p_1 \end{array} \right) +\dots \\ && \qquad  \qquad \qquad 
= \left( \begin{array}{c}  \frac{\gamma (- \epsilon)^{N-1} p_1}{\omega^N}  \\ ( - \frac{\epsilon}{\omega} )^{N-1} p_1 \end{array} \right) + \dots
\end{eqnarray}
Using this result for $q_N^*$, we construct the leading order terms for $q_j^*$ for $j=N-1, N-2, \dots , 2$.   For instance, consider the equation
\begin{eqnarray}
    0&=&f_{N-1}= \epsilon (q^*_{N-2}+q^*_N-2q^*_{N-1}) +
    \\ \nonumber && \qquad \qquad \qquad 
    + \omega q_{N-1}^* +((p_{N-1}^*)^2 +(q_{N-1}^*)^2)q^*_{N-1}
\end{eqnarray}
Inserting the leading order expression for $q_N^*$
from \eqref{eq:q*p*}, we find
\begin{eqnarray}
(-1)^{N-1} \gamma \left( \frac{\epsilon}{\omega}\right)^N p_1 + (\omega -2\epsilon) q_{N-1}^* +{\cal{O}}(\epsilon^{N+1})= 0\ ,
\end{eqnarray}
or
\begin{equation}
    q_{N-1}^* = (-1)^N \frac{\gamma}{\omega} \left( \frac{\epsilon}{\omega}\right)^N p_1
    +{\cal{O}}(\epsilon^{N+1})\ .
\end{equation}
Continuing in this way, we find the leading order expression
\begin{equation}\label{eq:qj}
q_j^* = (-1)^{2N-j-1} \frac{\gamma}{\omega} \left( \frac{\epsilon}{\omega} \right)^{2N-j-1} p_1 + \dots
\end{equation}
for $j=2, \dots , N-1$.

If we now substitute the leading
order expression for $q_2^*$ into the equation $f_{N+1} = 0$, (see \eqref{eq:beta} in Appendix \ref{ap:IFT} for the form of $f_{N+1}$)
  we calculate that the unique value of $\beta$ up to leading order is:
\begin{equation}
    \beta = (-1)^{2N-2} \gamma \left( \frac{\epsilon}{\omega} \right)^{2N-2} + \dots
    = \gamma \left( \frac{\epsilon}{\omega} \right)^{2N-2} + \dots
\end{equation}

\begin{remark} The last equality used the fact that $2N-2$ is always even and it shows that
the value of $\beta$ is positive for any choice of $N$, and also extremely small when $N$ becomes large.  The positivity of $\beta$ is important and natural, since in order to have a steady state in the damped system we expect to need to inject energy at some point, which corresponds to a positive value of $\beta$.
\end{remark}

\subsection{The ``twist'' of the breathers}\label{sec:twist}

For the undamped breathers, the fixed points all have $q_j^*=0$ which means that
the $z_j$ all lie on the real line within the complex plane - although by
rotational invariance, they can be rotated to any line through the origin.
In contrast, the damped and driven breathers calculated in the previous section have
$q_j^* \ne 0$ for $j>1$, when $q_1=0$. Therefore, the  components of the damped and driven
breathers do not lie on the same line through the origin. This is the twist
of the breathers which we will explore in this section.

An appropriate form of analyzing the twist is by studying the norm and the
complex phase of the equations $z_j(t) = r_j(t) e^{i\varphi_j(t)}$.
The utility of such a representation
in the study of periodic solution of the undamped dNLS equation was
noted at least as long ago as \cite{Eilbeck:1985}.  Since we are interested in studying the energy of the system, we can look
at the energy of each site individually, defined by $E_j(t) = \tfrac{1}{2} |z_j(t)|^2$.
The system can be represented using these energies and the complex phase.
Furthermore, due to rotational invariance we can use the differences in the phases rather
than the absolute phases to represent the system. The breather solutions
correspond to stationary points of this representation. Calling
$\psi_j(t) = \varphi_{j+1}(t) - \varphi_j(t)$ for $j=1, 2, \dots, N-1$, the system can be written as follows:
  
\begin{eqnarray}\label{eq:Ej}
\dot{E}_j &=& - 2 \epsilon \sqrt{E_j E_{j-1}} \sin \psi_{j-1} + 2 \epsilon \sqrt{E_{j+1} E_j} \sin \psi_j +2 \beta \delta_{1,j} E_1 - 2 \gamma \delta_{N,j} E_N\\
\dot{\psi}_1 &=& 2 ( E_2 - E_1 ) \left( 1 + \epsilon \frac{\cos\psi_1}{2 \sqrt{ E_1 E_2}} \right) - \epsilon \sqrt{\frac{E_3}{E_2}} \cos \psi_2 + \epsilon\\
\dot{\psi}_j &=& 2 ( E_{j+1} - E_j ) \left( 1 + \epsilon \frac{\cos \psi_j}{2 \sqrt{E_j E_{j+1}}} \right) - \epsilon \left( \sqrt{\frac{E_{j+2}}{E_{j+1}}} \cos \psi_{j+1} - \sqrt{\frac{E_{j-1}}{E_j}} \cos \psi_{j-1} \right)\\
\dot{\psi}_{N-1} &=& 2 ( E_N - E_{N-1} ) \left( 1 + \epsilon \frac{\cos \psi_{N-1}}{2 \sqrt{E_{N-1} E_N}}\right) + \epsilon \sqrt{\frac{E_{N-2}}{E_{N-1}}} \cos \psi_{N-2} - \epsilon
\end{eqnarray}
 
\begin{remark} Note that an additional advantage of the energy-phase representation is that it uses the rotational invariance of
the system to reduce the dimension of the system of equations from
$2N$ to $2N-1$.
\end{remark}
\begin{remark} In \eqref{eq:Ej}, we set the first term equal to zero if $j=1$, and the second term equal to zero if $j=N$.
\end{remark}

This energy-phase representation of the breathers is especially
useful when we wish to compare the metastable states which emerge 
when the damped system is allowed to decay with our breather solutions.
In particular, 
in the undamped,  n-site case all the breathers have $\sin \psi_j = 0$ for all $j$.  This means, for example, that if we use the rotational invariance
to set $q_1^*=0$, then all the $q_j^*$ will be equal to 
zero.  However, for nonzero
$\beta,\gamma$ the equations will lead to $\sin \psi_j \ne 0$. 

For the particular breather
calculated in Section \ref{sec:approx}, the phase is calculated from $\varphi_j = \arctan ( \tfrac{q_j^*}{p_j^*} ) = \tfrac{q_j^*}{p_j^*} + \dots$.
The phases for this breather will be:
\begin{equation}
    \varphi_j = \frac{\gamma}{\omega} \left( - \frac{\epsilon}{\omega} \right)^{2 ( N - j )} + \dots = \frac{\gamma}{\omega} \left( \frac{\epsilon}{\omega} \right)^{2(N-j)} + \dots\ , {\mathrm{for}}\ \ j= 2, 3 \dots , N\ ,
\end{equation}
while $\varphi_1 = 0$ as a consequence of choosing $q_1=0$. Using this, we find that $\psi_j = \varphi_{j+1} - \varphi_j$ to leading order is simply:
\begin{equation}
    \psi_j = \frac{\gamma}{\omega} \left( \frac{\epsilon}{\omega} \right)^{2(N-j-1)} + \dots
\end{equation}

In particular, we obtain the interesting following result for the last phase difference:

\begin{theorem}\label{th:twist_approx}
For a fixed $\gamma>0,\omega>0$, the phase difference between the last two oscillators, $\psi_{N-1}$, is up to leading order independent of the number of sites. The expression for this phase is:
\begin{equation}\label{eq:twist_prediction}
    \psi_{N-1}=\frac{\gamma}{\omega}+\dots
\end{equation}
\end{theorem}

\begin{remark} As we will see in the numerical experiments in Section \ref{sec:numerics}, the metastable states we observe during the energy decay in these systems, also exhibit a non-zero ``twist'' of this type, whose value is very close to the prediction of \eqref{eq:twist_prediction}.
\end{remark}

\section{Stability} \label{sec:spec} 

In this section, we examine the stability of the damped and driven breathers constructed above, using the fact that they are fixed points of the equations \eqref{eq:NLS_P}-\eqref{eq:NLS_Q}.  Define
\begin{eqnarray}\label{eq:FP}
&& \tilde{F}(p_1, \dots , p_N, q_1, \dots , q_N, \beta ; \epsilon, \omega, \gamma) \\ \nonumber && \qquad \qquad 
= (f_1, f_2, \dots , f_N, f_{N+1}, \dots f_{2N})\ ,
\end{eqnarray}
where the component functions $f_j$ are defined in \eqref{eq:f2-Nm}-\eqref{eq:fN}-\eqref{eq:f1} in Appendix \ref{ap:IFT}.
We recall that our goal in this stability analysis is to show that the stability properties
of the family of breathers constructed in the previous section is consistent with the trapping
of trajectories observed in studies of heat conduction.

If we linearize $\tilde{F}$ at the the fixed points $(p^*,q^*)=(p_1^*, \dots , p_N^*,q_1^*,\dots , q_N^*)$, 
the Jacobian matrix takes the form
\begin{equation}
D_{(p,q)} \tilde{F}\big|_{(p^*,q^*)} = \left( \begin{array}{cc} D^{(1)} & A^{(\omega,\epsilon)} \\ B^{(\omega, \epsilon)} &  D^{(2)} 
\end{array} \right)
\end{equation}

Here, $A^{(\omega,\epsilon)}$ and $B^{(\omega, \epsilon)}$ are $N\times N$ tri-diagonal matrices
and $D^{(1)}$ and $D^{(2)}$ are $N \times N$ diagonal matrices.  More precisely, 
the matrix $A^{(\omega,\epsilon)}$ has $\epsilon$ on the sub- and super-diagonal, while the
diagonal elements are $\omega -  a_j \epsilon -  (p_j^*)^2 - 3 (q_j^*)^3$, $j=1, 2, \dots , N$,
 where $a_1=1$, $a_N =1$ and all other $a_j = 2$.  The matrix $B^{(\omega, \epsilon)}$
 has sub- and super-diagonal elements equal to $-\epsilon$, while the diagonal
 elements are $-\omega + a_j \epsilon + 3(p_j^*)^2 + (q_1^*)^2$.  The diagonal matrix $D^{(1)}$ has
 diagonal elements $D^{(1)}_{jj} = \beta \delta_{1,j} - \gamma \delta_{N,j} -2 p_j^* q_j^*$, while
 $D^{(2)}$ has matrix elements $D^{(2)}_{jj} = \beta \delta_{1,j} - \gamma \delta_{N,j} +2 p_j^* q_j^*$.
 Note that from the approximations to $p_j^*$ and $	q_j^*$ derived in the preprint, 
 we see that $D^{(k)}_{jj} = - \gamma \delta_{j,N} + {\cal{O}}(\gamma \epsilon^{2(N-1)})$.  Thus,
 the diagonal part of the Jacobian is dominated by the dissipation coming from $\gamma$.
 
 We begin by recalling  a few key facts that we already know about the spectrum of the linearization
 which we will need in order to approximate the ``small'' eigenvalue.  The proofs of these facts can be found in \cite{Eckmann:2020}
 \begin{enumerate}
 \item Due to the invariance of the original system of equations under complex rotations, there
 is a zero eigenvalue with eigenvector
 \begin{equation}
\vv^{(1)} = \left( \begin{array}{c} -q^* \\ p^* \end{array} \right)\ .
\end{equation}
\item The linearization also has $2N-2$ {\em simple} eigenvalues near $\pm i$ which all have negative real parts that are $\OO(\gamma)$.  Denote these eigenvalues by $\lambda_3, \lambda_4, \dots , \lambda_{2N}$ and the associated eigenvectors by $\vv^{(j)}$, $j = 3, 4, \dots , 2N$.  Note
that in \cite{Eckmann:2020} this is proven for the linearization about the breather with $\beta = 0$.  However, since all of these eigenvalues
are simple, this remains true for $\beta$ sufficiently small by simple, non-degenerate perturbation theory estimates.
\end{enumerate}

We also identify another vector 
\begin{equation}
\tvv^{(2)} = \left( \begin{array}{c} \partial_{\omega} p^* \\ \partial_{\omega} q^* \end{array} \right)\ ,
\end{equation}
which corresponds to the tangent vector in the direction of changes in frequency of the breathers.
While not an eigenvector, there are two important properties of $\tvv^{(2)}$ that we use in what follows:
\begin{enumerate}
\item The set of vectors $\{ \vv^{(1)}, \tvv^{(2)}, \vv^{(3)}, \vv^{(4)}, \dots , \vv^{(2N)} \}$ form
a basis for $\real^{2N}$.  This follows from the fact that they are small perturbations of the eigenvectors
and generalized eigenvector of the linearization of the undamped breathers studied in \cite{Eckmann:2020}.
\item We know that when the linearization at the breather acts on $\tvv^{(2)}$ we have
\begin{equation}\label{eq:tv2action}
D_{(p,q)} \tilde{F}\big|_{(p^*,q^*)} \tvv^{(2)} = \vv^{(1)} - (\partial_{\omega} \beta) \left( \begin{array}{c} p_1^* \\
0 \\ \cdot \\ \cdot \\ \cdot \\ 0 \end{array} \right)\ .
\end{equation}
This follows by differentiating the equation for the damped and driven breather with respect to $\omega$.
\end{enumerate}

From our approximations to the breather we have
\begin{equation}
\partial_{\omega} p_1^* = \frac{1}{2\sqrt{\omega}}+ \dots\ , \ {\mathrm{and}}\  \ \partial_{\omega} \beta = 
- (2N-2) \frac{\gamma}{\omega} \left(\frac{\epsilon}{\omega} \right)^{2N-1} + \dots
\end{equation}

From the form of $\tvv^{(2)}$, and the expression for $\partial_{\omega} p_1^*$ we have
\begin{equation}
\left( \begin{array}{c} p_1^* \\
0 \\ \cdot \\ \cdot \\ \cdot \\ 0 \end{array} \right) = 2 \omega \tvv^{(2)} + {\OO(\epsilon)}\ .
\end{equation}

\begin{remark}  To simplify the following expressions, we choose $\omega =1$  for the remainder
of this section, but the case of general $\omega$ is very similar.
\end{remark}

Then, if we expand $(p_1^*,0,\dots,0)^T$ in terms of the basis
$\{ \vv^{(1)}, \tvv^{(2)}, \vv^{(3)}, \vv^{(4)}, \dots , \vv^{(2N)} \}$,  we have
\begin{equation}
\left( \begin{array}{c} p_1^* \\
0 \\ \cdot \\ \cdot \\ \cdot \\ 0 \end{array} \right) = \alpha_2 \tvv^{(2)} + \alpha_1 \vv^{(1)}
+ \sum_{j=3}^{2N} \alpha_j \vv^{(j)}\ ,
\end{equation}
with $\alpha_2 = 2 + \dots$, and all the other $\alpha_j$ at least $\OO(\epsilon)$ or smaller.
Note that the coefficients $\alpha_j$ are, in principle,  computable, at least approximately, since we have explicit
formulas for $\vv^{(1)}$ and $\tvv^{(2)}$, and since the remaining $\vv^{(j)}$'s correspond to simple
eigenvalues, they can be computed approximately by ordinary non-degenerate perturbation theory, at
least to leading order.

We now denote the small, non-zero eigenvalue as $\lambda_2$, and write
its eigenvector as
\begin{equation}\label{eq:v2exp}
\vv^{(2)} = \vv^{(1)} + \mu_2 \tvv^{(2)} + \sum_{j=3}^{2N} \mu_j \vv^{(j)}\ .
\end{equation}
Note that we compute below that the $\mu_j$'s are all small and go to zero as $\gamma \to 0$.
This is expected since for $\gamma =  0$, we know that the linearization at the breather 
has a two-dimensional zero eigenspace, but only one zero eigenvector, $\vv^{(1)}$ \cite{Eckmann:2020}.

Now consider the eigenvalue equation for $\lambda_2$:
\begin{equation}\label{eq:eig2}
D_{(p,q)} \tilde{F}\big|_{(p^*,q^*)} \vv^{(2)} = \lambda_2 \vv^{(2)}.
\end{equation}
Insert the expansion from \eqref{eq:v2exp} into both sides of this equation and use the fact that
on the left hand side we know that
$D_{(p,q)} \tilde{F}\big|_{(p^*,q^*)} \vv^{(1)} = 0$, $D_{(p,q)} \tilde{F}\big|_{(p^*,q^*)} \vv^{(j)}
= \lambda_j \vv^{(j)}$, for $j = 3,  \dots , 2N$, and 
$D_{(p,q)} \tilde{F}\big|_{(p^*,q^*)} \tvv^{(2)}$ is given by \eqref{eq:tv2action}.

Thus, we have
\begin{eqnarray}\nonumber
&& \mu_2 \left(\vv^{(1)} - (\partial_{\omega} \beta) \left( \alpha_2 \tvv^{(2)} + \alpha_1 \vv^{(1)}
+ \sum_{j=3}^{2N} \alpha_j \vv^{(j)} \right)  \right) + \sum_{j=3}^{2N} \lambda_j \mu_j \vv^{(j)}  = \\ \label{eq:eig2expanded}
&& \qquad \qquad \qquad \qquad  = \lambda_2 \left( \vv^{(1)} + \mu_2 \tvv^{(2)} + \sum_{j=3}^{2N} \mu_j \vv^{(j)} \right)\ .
\end{eqnarray}
We now equate the coefficients of the various basis vectors on the two sides of this equation.
From the coefficients of the vector $\tvv^{(2)}$, we have
\begin{equation}\label{eq:v2coeff}
-\alpha_2 \mu_2 \partial_{\omega} \beta = \lambda_2 \mu_2\ ,\ \ {\mathrm{or}}\ \ \lambda_2 = 
-\alpha_2 (\partial_{\omega} \beta) \ .
\end{equation}
Inserting our approximate expressions for $\alpha_2$ and $(\partial_{\omega} \beta)$, we find the approximate
expression for the small, positive eigenvalue:
\begin{equation}\label{eq:lambda2approx}
\lambda_2 = 2 (2N-2) \gamma \epsilon^{(2N-2)} + \dots
\end{equation}
(Recall that we have set $\omega=1$.)  While they are not really important, we can also solve for the remaining coefficients $\mu_j$.  Thus, equating the coefficients of $\vv^{(1)}$ in \eqref{eq:v2coeff},
we find
\begin{equation}
\mu_2 = \lambda_2/(1-\alpha_1 (\partial_{\omega} \beta))\ ,
\end{equation}
while
\begin{equation}
\mu_j = \frac{-\alpha_j \mu_2 (\partial_{\omega} \beta)}{\lambda_j - \lambda_2}\ ,\ \ j = 3, \dots , 2N\ .
\end{equation}

\begin{remark}  Note that for $N=2$, the formula in \eqref{eq:lambda2approx} agrees with the 
computation for the two-site breather in energy-angle coordinates in {\eqref{eq:l22s}} \end{remark}

From this expression, we obtain $\mu_j<<\mu_2$ for $j\geq 3$. Thus, from \eqref{eq:v2exp}, the small positive eigenvector is (up to very high order) in the  span of the vectors  $\vv^{(1)},\tvv^{(2)}$ which corresponds to the tangent space of the cylinder of breathers, with $\vv^{(1)}$ corresponding to rotations around the cylinder and $\tvv^{(2)}$ to the direction along the family of breathers.  Hence, the positive eigenvalue pushes solutions to slide along the family of breathers, causing $\omega$ to decrease over time. Movement perpendicular to the cylinder is negligible until $\epsilon\sim\omega$, at which point $\vv^{(2)}$ has a significant component in the normal direction ($\mu_j\sim\mu_2$) forcing the solution to finally exit the neighborhood of the cylinder.

 \section{Effects on evolution/Metastable behavior}\label{sec:meta}
 
 We next consider how the spectral picture derived in the previous subsection can be
 used to model the metastable evolution in the damped (but not driven) system of oscillators.
 We will see that we reproduce the picture obtained from the modulation approach of \cite{Eckmann:2020},
 by considering the evolution near one of the fixed points represented by the damped and driven
 breathers.  However, as we discuss, these damped and driven breathers provide a more accurate
 approximation to the metastable state than the undamped breathers used as an approximation
 in the prior work.
 
 To simplify the exposition, we introduce the following notation.  We set $u=(p,q)^T$, where 
 $p$ and $q$ are $N$-vectors.  Likewise, the damped and driven breather will be denoted by
 $u^* = (p^*,q^*)$ - depending on the context, we may also include the dependence of the 
 breathers on the parameters - i.e. $u^* = u^*(\gamma,\omega)$.  We also introduce two
 diagonal matrices which represent respectively the driving and damping - 
 $B(\beta)$ is a $2N\times 2N$ matrix with $(1,1)$ and $((N+1),(N+1))$ element equal to $\beta$
 and all other elements zero, and $\Gamma(\gamma)$ is $2N\times 2N$ matrix with $(N,N)$ and $(2N,2N)$ element equal to $\gamma$
 and all other elements zero.  Finally, we let $H(u)$ denote the Hamiltonian of the undamped and undriven
 lattice of oscillators.  We can then write the damped, but not driven, evolution as
 \begin{equation}\label{eq:damped_not_driven}
\dot{u} = J \nabla_u H(u) - \Gamma u\ .
\end{equation}
The damped and driven breathers $u^*$ are solution of the system of equations
\begin{equation}\label{eq:DDnotation}
J \nabla_u H(u^*) + B(\beta^*) u^* - \Gamma(\gamma) u^* =  0\ .
\end{equation}
Here, $J$ is the standard, symplectic matrix for Hamiltonian systems, which takes the form 
$J = \left( {\begin{array}{cc} 0 & 1 \\ -1 & 0 \end{array} }\right)$, with $0$ representing an $N\times N$ matrix of zeros, and $1$ an $N\times N$ 
identity matrix.

Recall that $u^*$ is the $\theta=0$ member of a family of fixed points of \eqref{eq:DDnotation} of
the form 
\begin{equation}\label{eq:theta_family}
U^*(\theta) =(\Re(e^{i \theta} u^*), \Im(e^{i\theta} u^*))\ .
\end{equation}
  Because of the equivariance
of the linearized equation under these complex rotations, the spectrum of the linearization
at the breather solution is independent of $\theta$, and hence, without loss of generality,
we can assume that $\theta=0$.

We will consider a solution of \eqref{eq:damped_not_driven} with initial conditions near the breather.
We will expand the solution with respect to the basis of vectors
$\{ \vv^{(1)}, \tilde{\vv}^{(2)}, \vv^{(3)} , \dots , \vv^{(2N)}\}$, where as in Section \ref{sec:spec},  $\vv^{(1)} = (-q^*,p^*)$
is the eigenvector of the linearization about the damped and driven breather with eigenvalue
zero, $\vv^{(j)}$, with $j=3, 4, \dots , 2N$ are the eigenvectors with eigenvalues near $\pm i$ (but
all with real parts $\sim -{\mathcal{O}}(\gamma)$) and $\tilde{\vv}^{(2)} = \partial_{\omega} u^*$.
We know from our discussion of the spectrum of the linearization about the damped and driven breather
that the eigenvector with small positive eigenvalue is approximately
$\vv^{(1)} + \mu_2 \tilde{\vv}^{(2)}$, with $\mu_2 \sim {\mathcal{O}}(\epsilon^{2N-2})$.

We take an initial condition for \eqref{eq:damped_not_driven} of the form
\begin{equation}\label{eq:decomp}
u_0 = u^* + w_0\ .
\end{equation}
We assume that $w_0$ is:
\begin{enumerate}
\item small, and
\item has no component in the $\vv^{(1)}$ or $\tvv^{(2)}$ direction.
\end{enumerate}

The second of these conditions can be enforced by adjusting the angle $\theta $ in \eqref{eq:theta_family}
and the frequency $\omega$.  Because of the equivariance of the linearization about the damped
and driven breather with respect to the complex rotations by $\theta$, the spectrum of the linearization
is independent of $\theta$ (and depends only weakly on $\omega$) so we will assume that we
are linearizing about the breather with $\theta=0$ and $\omega=1$.

We look for a decomposition of the solution similar to that in  \eqref{eq:decomp} - that is, we write the
solution of \eqref{eq:damped_not_driven} as
\begin{equation}
u = u^*+w\ .
\end{equation}
Inserting this into \eqref{eq:damped_not_driven} we find
\begin{eqnarray}\nonumber
\dot{w} &=& J \nabla H( u^* + w) - \Gamma (u^*+w) \\ \label{eq:weq1}
&=& J \nabla H(u^*) - \Gamma u^* + B u^* \\ \label{eq:weq2}
&& \qquad + \left(J \nabla H( u^* + w) - J \nabla H(u^*)\right) + B w - \Gamma w \\ \label{eq:weq3}
&& \qquad \qquad - B(u^* + w)
\end{eqnarray}
The goal in rewriting the equation for $\dot{w}$ in this form is to note that from the definition of
$u^*$, the expression in \eqref{eq:weq1} vanishes, while if we apply Taylor's Theorem
to expand \eqref{eq:weq2} in $w$, the linear terms in $w$ just give $\Lu w$, the linearization
about the damped and driven breather.
Thus, we have
\begin{equation}\label{eq:wdot}
\dot{w} = \Lu w - B w - B u^* + \tilde{N}(w)\ ,
\end{equation}
where $\tilde{N}(w)$ collects the nonlinear terms resulting from applying Taylor's Theorem to \eqref{eq:weq2},
and is ${\cal{O}}(w^2)$.

We can write the solution as
\begin{equation}
w(t) = \alpha_1(t) \vv^{(1)} + \talpha_2(t) \tvv^{(2)} + \sum_{j=3}^{2N} \alpha_j(t) \vv^{(j)}\ .
\end{equation}

We will ignore the contributions of $\alpha_j(t)$, with $j \ge 3$, since we expect the components of
the solution in those directions to be contracting since the eigenvalues corresponding to these
eigendirections all have negative real part.

To isolate the evolution of $\alpha_1(t)$ and $\talpha_2(t)$, we also introduce the vectors
\begin{equation}
\tnn^{(1)} = c_1 \left( \begin{array}{c} \partial_{\omega} q^* \\ \partial_{\omega} p^* \end{array} \right)
\ ,\ \ 
\nn^{(2)} = c_2 \left( \begin{array}{c} p^* \\ -q^* \end{array} \right)\ ,
\end{equation}
for constants $c_1$ and $c_2$ determined below.
Note that if we consider the transpose of $\Lu$, we have
\begin{equation}\label{eq:noneaction}
(\Lu)^T \tnn^{(1)} =  - \left( \begin{array}{c} 0 \\ \left( \begin{array}{c} (\partial_{\omega}\beta)  p_1^* \\ 0 \end{array} \right) \end{array} \right) + \frac{c_1}{c_2} \nn^{(2)} + \OO(\gamma \epsilon^{2N-1} )\ \ ,
\end{equation}
and
\begin{equation}
(\Lu)^T \nn^{(2)} = \left( \begin{array}{c}  0 \\ 0 \end{array} \right) + \OO(\gamma \epsilon^{2N-1} )\ .
\end{equation}

The first vector on the right hand side of \eqref{eq:noneaction} is the $2N$ vector whose
first $N$ components are all zero, and whose second $N$ components have
$(\partial_{\omega}\beta)  p_1^*$ in the $(N+1)^{st}$ component and then zero in the remaining
components.
We choose the constants $c_1$ and $c_2$ so that
$\tnn^{(1)} \cdot \vv^{(1)} =\nn^{(2)} \cdot \tvv^{(2)} = 1$.  Note that from the definitions of these
vectors, and the asymptotic formulas we derived for $p^*$ and $q^*$, we find that
for $\omega=1$, $c_1 = 2 + \OO(\epsilon)$, $c_2=2+\OO(\epsilon)$ and that
$\tnn^{(1)} \cdot \tvv^{(2)}$ and $\nn^{(2)} \cdot \vv^{(1)}$ are both $\OO(\gamma \epsilon^{2N-1})$.

Thus, if we take the dot product of the equation \eqref{eq:wdot} with $\tnn^{(1)}$, we have (ignoring
the nonlinear terms, for the moment)
\begin{eqnarray}
\dot{\alpha}_1 &=& \alpha_1 \langle \tnn^{(1)}, \Lu \vv^{(1)} \rangle 
+ \tilde{\alpha}_2 \langle \tnn^{(1)}, \Lu \tvv^{(2)} \rangle 
\\ \nonumber
&& \qquad \qquad - \langle \tnn^{(1)},B u^* \rangle 
- \langle \tnn^{(1)},B (\alpha_1 \vv^{(1)} + \tilde{\alpha}_2 \tvv^{(2)} ) \rangle
- \dot{\tilde{\alpha}}_2 \langle \tnn^{(1)} , \tvv^{(2)} \rangle \\ \nonumber
&=& \alpha_1 \langle \left( \begin{array}{c} 0 \\ \left( \begin{array}{c} (\partial_{\omega}\beta)  
p_1^* \\ 0 \end{array} \right) \end{array} \right) + \frac{c_1}{c_2} \nn^{(2)} , \vv^{(1)} \rangle 
+ \tilde{\alpha}_2 \langle  \left( \begin{array}{c} 0 \\ \left( \begin{array}{c} (\partial_{\omega}\beta)  
p_1^* \\ 0 \end{array} \right) \end{array} \right) + \frac{c_1}{c_2} \nn^{(2)},  \tvv^{(2)} \rangle 
\\ \nonumber
&& \qquad \qquad - \langle \tnn^{(1)},B u^* \rangle 
+ \OO(\epsilon^{2N-1} (\alpha_1+\tilde{\alpha}_2+\dot{\tilde{\alpha}}_2)) \\ \nonumber
&=&\frac{ c_1 \tilde{\alpha_2}}{c_2} + \OO(\gamma \epsilon^{2N-1} (\alpha_1+\tilde{\alpha}_2+ \dot{\tilde{\alpha}}_2))
+ \OO(\epsilon^{4N-2}) \ ,
\end{eqnarray}
or
\begin{equation}
\dot{\alpha}_1 =  (1+ \OO(\epsilon)) \tilde{\alpha}_2 + \OO(\gamma \epsilon^{2N-1} (\alpha_1+\tilde{\alpha}_2+\dot{\tilde{\alpha}}_2))
+ \OO(\epsilon^{4N-2}) \
\end{equation}
Here we used the form of $\tnn^{(1)}$ and $u^*$ to bound $| \langle \tnn^{(1)},B u^* \rangle |$
by $\OO(\epsilon^{4N-2})$.

If we now repeat this procedure, taking the dot product of \eqref{eq:wdot} with $\nn^{(2)}$, we find
\begin{eqnarray}
\dot{\tilde{\alpha}}_2 &=& \alpha_1 \langle \nn^{(2)}, \Lu \vv^{(1)} \rangle 
+ \tilde{\alpha}_2 \langle \nn^{(2)}, \Lu \tvv^{(2)} \rangle  - \dot{\alpha}_1 \langle \nn^{(2)} , \vv^{(1)} \rangle
\\ \nonumber
&& \qquad \qquad - \langle \nn^{(2)},B u^* \rangle 
- \langle \nn^{(2)},B (\alpha_1 \vv^{(1)} + \tilde{\alpha}_2 \tvv^{(2)} ) \rangle \\ \nonumber
&=&  - \langle \nn^{(2)},B u^* \rangle 
 + \OO(\gamma \epsilon^{2N-1} (\alpha_1+\tilde{\alpha}_2+\dot{\alpha}_1)) \\ \nonumber
 &=& -c_2 \beta + \OO(\gamma \epsilon^{2N-1} (\alpha_1+\tilde{\alpha}_2 +\dot{\alpha}_1))\ , \\ \nonumber
\end{eqnarray}
or
\begin{equation}
\dot{\tilde{\alpha}}_2 = -(2 +\OO(\epsilon)) \beta + \OO(\gamma \epsilon^{2N-1} (\alpha_1+\tilde{\alpha}_2+\dot{\alpha}_1))\ .
\end{equation}

If we ignore the high order correction terms, we can solve these equations and we find
\begin{equation}
\tilde{\alpha}_2(t) = \tilde{\alpha}_2(0) -2 \beta t\ ,
\end{equation}
and
\begin{equation}
\alpha_1(t) = \alpha_1(0)+\tilde{\alpha}_2(0) t - \beta t^2\ .
\end{equation}
By adjusting the breather $u^*$ about which we are perturbing, we can insure that $\alpha_1(0)
= \tilde{\alpha}_2(0) =0$, so we have
\begin{equation}
\alpha_1(t) = - \beta t^2\ ,\ \ {\mathrm{and}}
\ \  \tilde{\alpha}_2(t) = -2 \beta t\ .
\end{equation}
Inserting the value of $\beta$, we have
\begin{equation}
\alpha_1(t) = - \gamma \epsilon^{2N-2}  t^2\ ,\ \ {\mathrm{and}}
\ \  \tilde{\alpha}_2(t) = -2 \gamma \epsilon^{2N-2} t\ .
\end{equation}
Recall that our parameter $\gamma$ corresponds to $\gamma \epsilon$ 
and our $\tilde{\alpha}_2$ corresponds to $\phi$ in the paper
\cite{Eckmann:2018}.  (This last equivalence follows from the fact that $\phi$ in \cite{Eckmann:2018}
is the change in the frequency of the breather, while our $\tilde{\alpha}_2$ corresponds to the
displacement in the direction $\partial_{\omega} u^*$ - i.e. also the change in frequency of the breather.)
Then for $N=3$, the case considered in that reference, the formula for $\phi(t)$ derived 
in that work agrees exactly with our formula for $\tilde{\alpha}_2$. 

To compare the expression for $\alpha_1(t)$ with that of $\theta(t)$ in \cite{Eckmann:2018}, note
that $\alpha_1(t)$ is the coefficient of $\vv^{(1)}$, the eigenvector corresponding to 
rotations in the complex plane $u \to e^{i \theta} u$.  (This is not the same $\theta$ as in 
\cite{Eckmann:2018}).  In \cite{Eckmann:2018}, the
quantity corresponding to a complex rotation of this sort is $t \phi(t) + \theta(t)$  (see equation (13) of that reference).
If we take the formulas $\theta(t)$ and $\phi(t)$ from \cite{Eckmann:2018} and compute,
we obtain
\begin{equation}
t \phi(t) + \theta(t) = - 2 \gamma \epsilon^5 t^2 + \gamma \epsilon^5 t^2 =  -\gamma \epsilon^5 t^2\ ,
\end{equation}
which again, agrees exactly with our formula for $\alpha_1(t)$ when $N=3$.  (Again, recalling
that our $\gamma$ should be replaced by $\gamma \epsilon$ to compare with
the results of \cite{Eckmann:2018}.)

This is only an approximate calculation - in particular, we have ignored the contributions of the 
nonlinear terms, but in \cite{Eckmann:2020} it was shown rigorously that in the case of damping alone, the nonlinear effects could be
rigorously incorporated without changing the predictions of the linear approximation and we expect that the same is true here.

\section{Numerical verification}  \label{sec:numerics}

In this section, we present numerical experiments which illustrate
the results of the preceding sections.  We will compute the spectrum of the linearization
about the breather, showing that it has the form
described in the previous section, and we will compute the dynamics of the 
 metastable state of the damped system.  In addition to showing that
 solutions of the {\em damped} system, which start close to the damped and driven
 breathers remain close to this family for a long time we will also follow such
 solutions until their disappearance and develop a conjecture about how such
 families of metastable states terminate.
We focus on systems with a small number
of sites (specifically, $N=2,3$) because the very slow drift along the family
of breathers means that the computational time required increases very 
rapidly with $N$. {This makes it unfeasible
to compute the solution up until the disappearance of the
metastable state for larger $N$}.  In addition, by
concentrating on these small systems
we are able to perform many computations explicitly, which helps
to elucidate the numerical results.  All numerical computations are done
with Mathematica.

We begin by computing the spectrum of the linearization of equations \eqref{eq:NLS_P}-\eqref{eq:NLS_Q} about the breather
solution.  Table \ref{tab:spec} contains the spectrum of the linearization around the breather solution
with $\epsilon = 0.03$, $\gamma = 0.0035$, and $\omega =1$.  Note that as expected, there are two
eigenvalues near $ i$, and two more near $-i$, all with negative real part,   one with
small, positive real part and one eigenvalue which is zero to within the errors of our computations (i.e. $ \sim {\OO(10^{-9}})$),

\begin{table}\label{tab:spec}
\begin{center}
\caption{Numerically computed spectrum of the linearization about the breather, when $\epsilon = 0.03$, $\gamma = 0.0035$, and $\omega =1$.}
\begin{tabular}{c }
\hline
$-0.00257707 + 0.987135 i$ \\ $-0.00257707 - 0.987135 i $\\ $ -0.000922942 + 
 0.914551 i $\\ $-0.000922942 - 0.914551 i$ \\ $3.43099*10^{-8} + 0. i $\\
$2.28602*10^{-9} + 0. i $\\
\hline
\end{tabular}
\end{center}
\end{table}

\begin{figure}
    \centering
    \includegraphics[width=8cm]{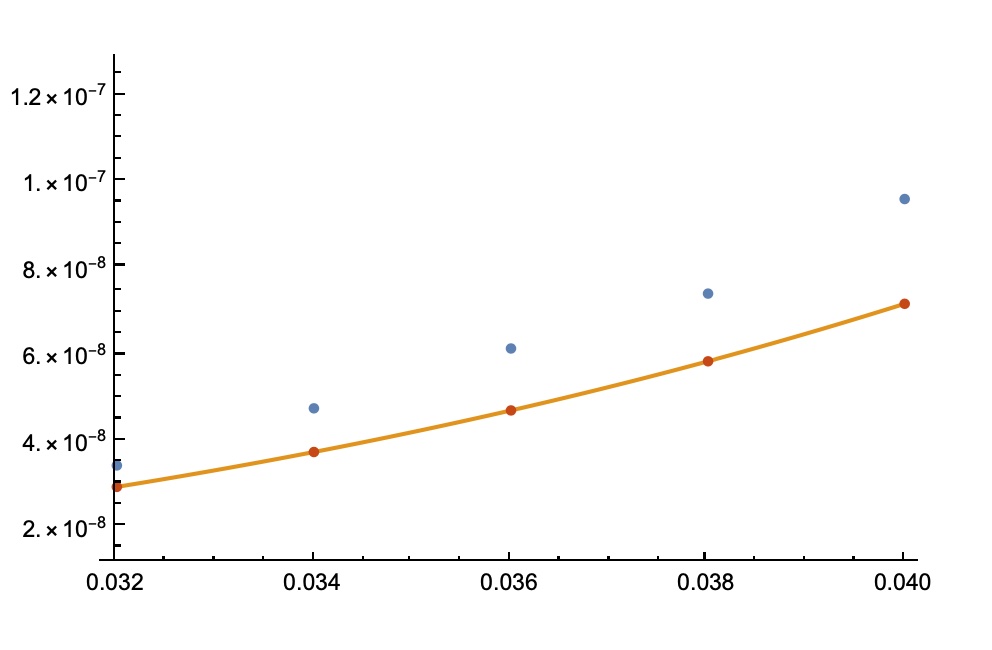}
    \caption{Numerically computed value of the small, positive eigenvalue (blue dots), when $N=3$, $\gamma=0.0035$, and $\epsilon_k=
    0.03+0.002*k$, for $k=1, 2, \dots, 5$, compared with the formula in \eqref{eq:lambda2approx}}
    \label{fig:small_eigenvalue}
\end{figure}

As we saw in the previous section, the crucial eigenvalue for understanding the metastable
behavior of the system is the small positive eigenvalue, so we focus on that.  In Figure \ref{fig:small_eigenvalue},
we plot the numerically computed value of this eigenvalue for $N=3$,  $\gamma = 0.0035$, and for five
different values of $\epsilon$, $\epsilon_k = 0.03+0.002*k$, $k=1, 2, \dots , 5$, and compare these values to 
the values predicted by equation \eqref{eq:lambda2approx}.  As we see, the computed values are quite close to the theoretically predicted values,
and the approximation improves as $\epsilon$ becomes smaller.

We next compute the full time evolution of the three site system to examine how long the metastable
behavior persists, and how it terminates.  For these computations we found it easier to work in the energy-phase
representation of the system,  \eqref{eq:Ej}, since this reduces the number of equations we need to solve by one.
We choose as initial conditions for our numerics the approximation to the 
damped and driven breather computed in Section \ref{sec:approx}.   
 We also need some way of determining
what member of the breather family we are close to at any given time in the course of the the
computation.     

We proceed as follows.  Choose a value of $\omega_0$ to determine the initial values for the
equations.   Then set $p_1 = \sqrt{\omega_0}$.
and approximate the initial values for the other variables using
the equations of motion and our knowledge of the breather solution.

As the solution evolves, it moves along the family of breathers, and we need to
determine how the frequency (which determines which breather we are close to) changes.
If we are very close to a damped and driven breather we can consider
$\dot{\psi}_j \approx 0$ and so all the $\dot{\varphi}_j$ are the same.
For this reason, we can define $\omega(t)$ to be any of the $\dot{\varphi}_j$
and we should obtain approximately the same result. For simplicity, we define
it as $\omega(t) = \dot{\varphi}_1(t)$ and in turn we can calculate it from:

\begin{equation}\label{eq: omega}
    \dot{\varphi}_1 = 2 E_1 + \epsilon - \epsilon \sqrt{\frac{E_2}{E_1}} \cos \psi_1
\end{equation}

Notice that $\omega_0 \ne \omega(0)$ because of the fact that we calculate an
approximation to $p_1$ from $\omega_0$ and then use it to calculate
approximations for $E_1(0), E_2(0), \psi_1(0)$. In fact $\omega_0 = 2E_1(0)$.

\begin{figure}
    \centering
    \includegraphics[width=6cm]{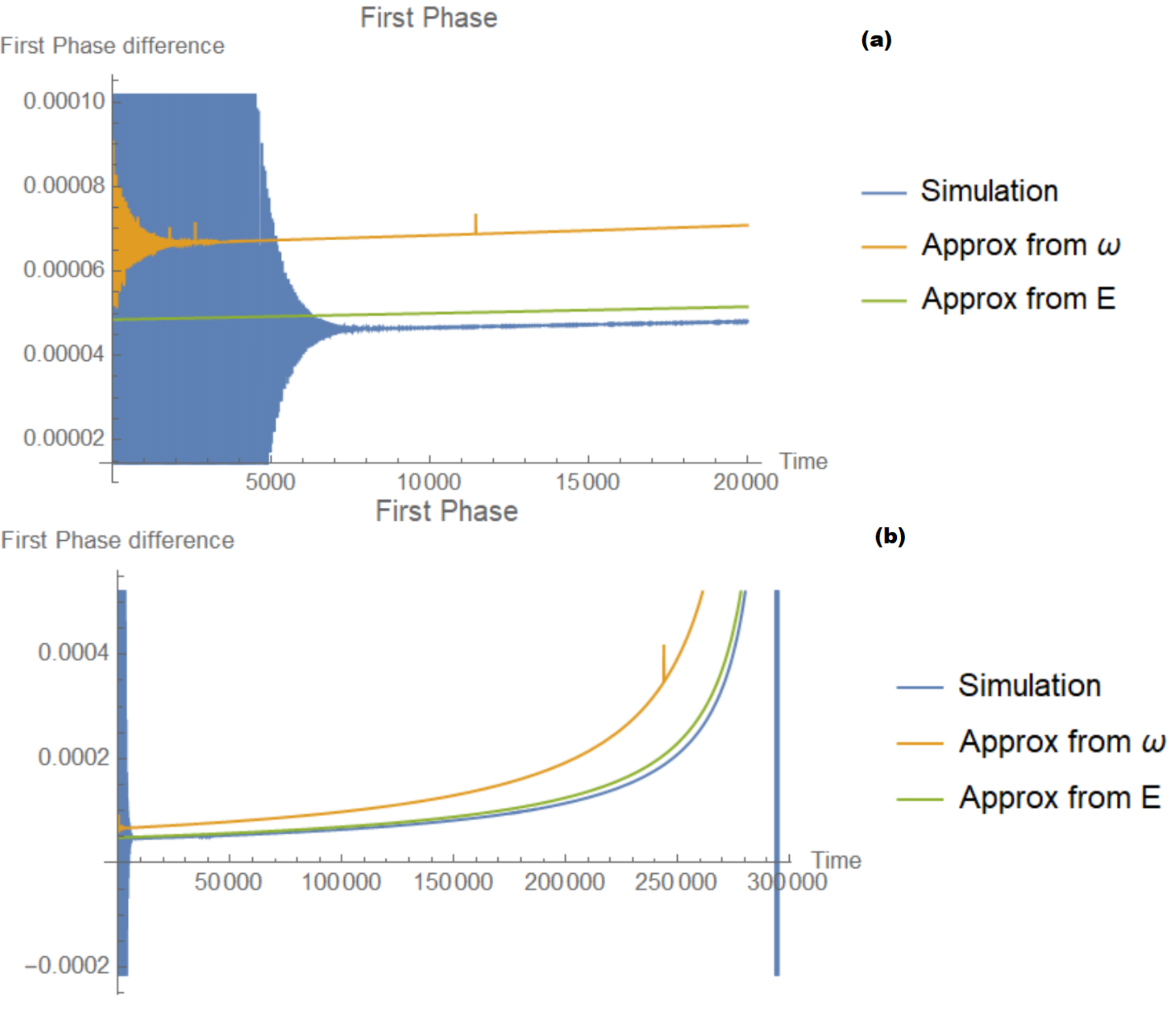}
    \caption{Graph of the phase difference between the first and second oscillator and its approximation for $\omega_0=1,\gamma=0.005,\epsilon=-0.1$ and initial condition given by the approximation to the damped and driven breathers of Section \ref{sec:approx}; (a) is graphed between times $0$ and $20\,000$, and (b) between times $0$ and $295\,000$, approximately when it leaves the metastable state.}
    \label{fig:3site phase1}
\end{figure}

The first quantity we examine are the graphs of $\psi_1 (t)$ and $\psi_2 (t)$ for the
three site system.  We recall that for the undamped breathers, all of the $\psi_j(t)$ are zero, while the
damped and driven breathers have non-zero ``twist''.
Figure \ref{fig:3site phase1} graphs $\psi_1 (t)$ for the three-site system with two approximations. The blue curve
represents the computed value of $\psi_j(t)$.   From Theorem \ref{th:twist_approx}, we know that the lowest order approximation to
the twist of the breather is $\psi_j(t) \approx \tfrac{\gamma}{\omega(t)}$, and this is plotted in orange, using the approximation 
$\omega(t) \approx \dot{\varphi}_1(t)$.  We obtain a good approximation to the 
computed twist, showing that the damped and driven breathers provide a better approximation to the metastable states than the
undamped breathers.  Finally, recalling that $\omega(t) \approx 2 E_1(t)$ we have also plotted (in green) the quantity
$ \tfrac{\gamma}{E_1(t)}$.  Interestingly, this quantity seems to provide an even better approximation to the numerically computed
twist though at the moment we have no theoretical explanation of why this should be the case.

Speaking colloquially, our numerics show that the metastable state slides along the family of damped and driven breathers.  This is consistent with our computations of the spectrum of the linearization which showed that the span of the eigenspace of the zero eigenvector and the eigenvector with
small positive eigenvalue was very close to the tangent plane of the family of damped and driven breathers.
Moreover, this allows for the modeling of the dynamical evolution of the metastable state of the damped system from the damped and driven breathers (especially the time evolution of the phases) for long times. We believe that the period for which this approximation is accurate is connected to the bifurcations of the undamped system.

\begin{figure}
    \centering
    \includegraphics[width=8cm]{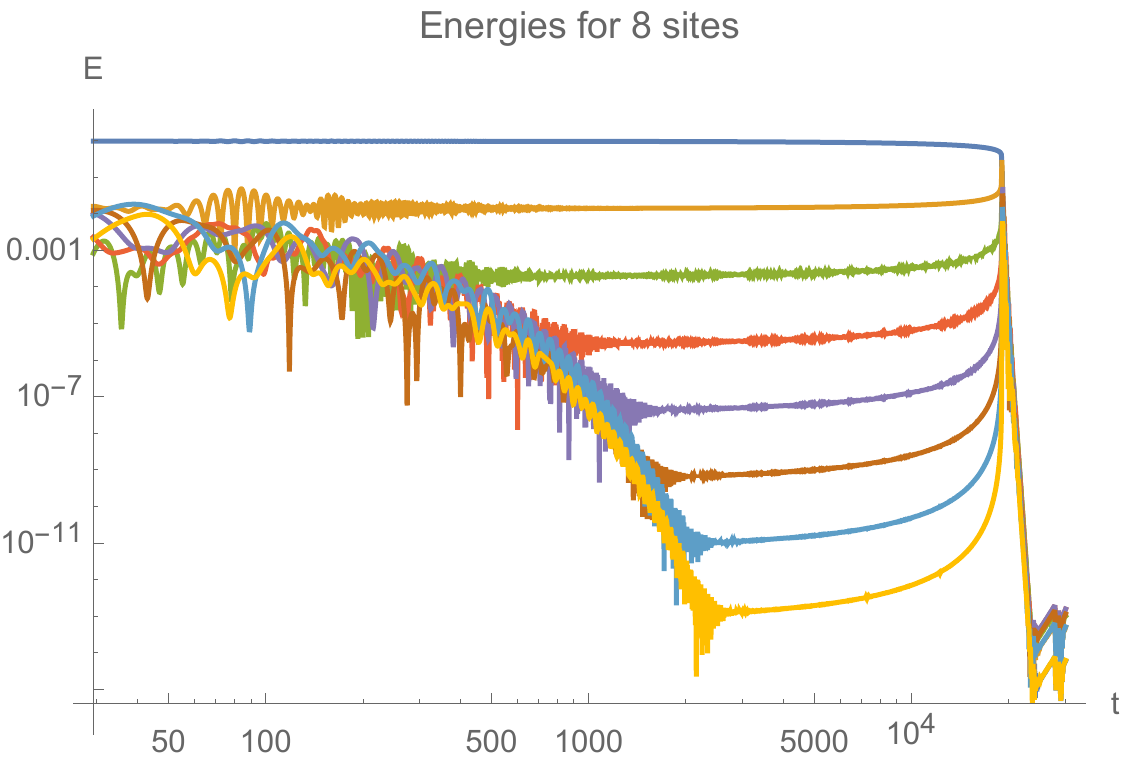}
    \caption{The energy of each oscillator in an eight-site system. The parameters chosen were
    $\gamma=0.05$, $\epsilon = -0.1$,  and the initial condition is $p_1=1$ and all other $p_j,q_j=0$.}
    \label{fig:8site energy}
\end{figure}

In Figure
\ref{fig:8site energy},  we present numerics  for a system with $N=8$ sites, which show that the effects we have described persist for
systems with larger numbers of oscillators.  Again we note that once the system has fully settled in the metastable state, we observe an initial, almost flat state where the system moves along the family of damped and driven breathers for a very long time. The system eventually reaches a point where the energy of the first oscillator decays suddenly while that of the other oscillators rise to meet it until they have the same 
energy, after which they reach the thermalized state and decay exponentially. The sudden decay in the energy is significantly more dramatic for the eight sites than it was for two or three sites.

Thus, we conjecture that the general evolution of these metastable states is as follows.  When the trajectory is ``captured'' by one of the damped and driven breathers, the system begins to evolve by drifting along this family until the family of breathers disappears through bifurcation, at which point the system may be captured by another attractive family of approximately periodic orbits.  There may, of course be other, possibly large, parts of the phase space where these metastable solutions do not exist and the evolution of the system is chaotic.  In particular, it is not clear what properties of the periodic orbit in the undamped Hamiltonian system result in the solution becoming attractive once the system is subjected to damping.

We now turn to a discussion of the two-site damped and driven system where many of the phenomena
discussed above can be illustrated by explicit computations.  A detailed derivation is found in Appendix \ref{ap: two-site analysis}, while in this section we will mostly show the results. The family of breathers of the Damped and Driven system related to the metastable state we have studied is analytically calculated to be:

\begin{equation}\label{eq:dd2}
    E_1=\frac12\sqrt{\frac\gamma\beta}\sqrt{\epsilon^2-\gamma\beta},\quad E_2=\frac12\sqrt{\frac\beta\gamma}\sqrt{\epsilon^2-\gamma\beta},\quad\psi = \Bigg\{
    \begin{array}{cc}- \pi + \arcsin\left( \frac{\sqrt{\gamma \beta}}{\epsilon} \right)&\quad  \epsilon >0\\
    -  \arcsin\left( \frac{\sqrt{\gamma \beta}}{\epsilon} \right) &\quad \epsilon <0.
\end{array}
\end{equation}

We can compute the spectrum of the linearization around these breathers, first writing out the Jacobian and then using Mathematica to derive series approximations. The small positive eigenvalue is:

\begin{eqnarray}\label{eq:l22s}
\lambda_2 &=& 4 \gamma \left( \frac{\epsilon}{\omega} \right)^2 + \mathcal{O}(\gamma\epsilon^4) \label{eq:eig1}
\end{eqnarray}

The zero eigenvalue, $\lambda_1$, is absent in the energy-phase representation, because
it corresponds to uniform phase rotation at all sites,  and the energy-phase representation considers only
phase differences between adjacent sites and hence is insensitive to phase rotations of the whole system.  The leading order
term in the small positive eigenvalue $\lambda_2$ has the same form as we derived in equation \eqref{eq:lambda2approx}, if
we set $N=2$ in that formula, as expected.

\begin{figure}
    \centering
    \includegraphics[width=7cm]{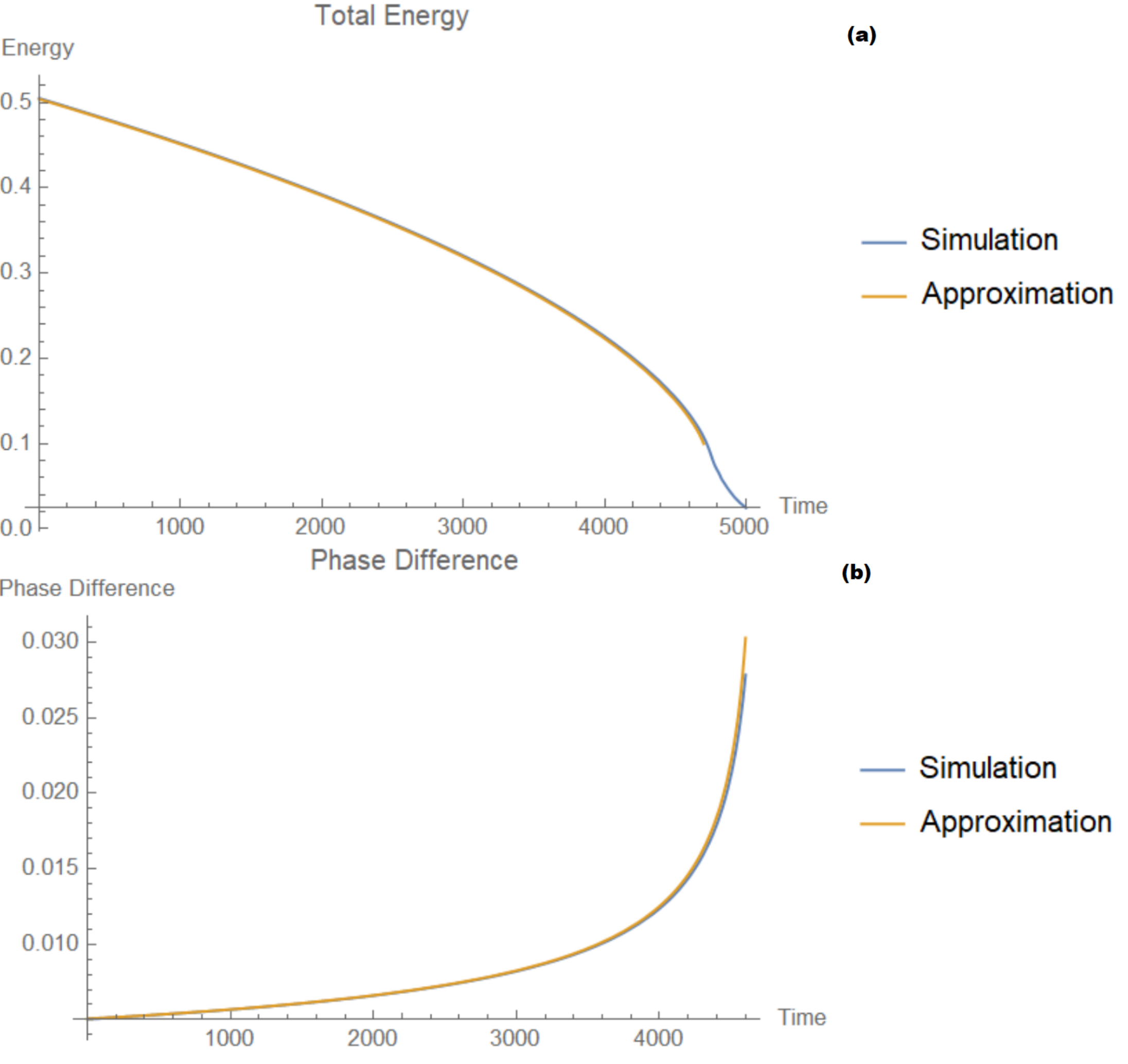}
    \caption{Graphs of the total energy and the phase difference for $E_0 = 0.505, \epsilon = -0.1, \gamma = 0.005$. (a) Is the total energy graphed from $t = 0$ to $t = \tau + 300 \approx 5000$ and (b) is the phase difference graphed from $t = 0$ to $\tau - 100 \approx 4600$.}
    \label{fig:2site approximation}
\end{figure}

In Appendix \ref{ap: two-site analysis}, we derive an approximation to the time evolution of the metastable state using the Ansatz that the solution maintains the relation $E_1(t)-E_2(t)=\sqrt{E(t)^2-\epsilon^2}$ resulting in an expression for the phase shift:

\begin{equation}\label{eq: 2siteshift}
    \sin \psi = \frac{\gamma}{2 \sqrt{E^2-\epsilon^2}}
\end{equation}
Note that this is only an approximation to the solution - this value
of  $\psi$ will not satisfy the equation for $\dot{\psi}$.
However Figure \ref{fig:2site approximation} shows that it is a very accurate approximation
when $E \gtrsim |\epsilon|$.

This result can then be used to calculate an approximation to the time evolution of the energy with initial energy $E(0)=E_0$:

\begin{equation}\label{eq: approximation}
    E(t) \approx \frac{|\epsilon|}{2} \left( \sqrt{-W_{-1} \left( - e^{-1 + 
    4\gamma (t - \tau)} \right)} + \frac{1}{\sqrt{-W_{-1} 
    \left( - e^{-1 + 4\gamma (t - \tau)} \right)}} \right)
\end{equation}

Here $W_{-1}$ is the $-1$ branch of the product logarithm, and 
$\tau$ is the time for which $\lim_{t \to \tau} E(t) = |\epsilon|$ and is given by:
\begin{equation}\label{eq:tau}
    \tau = \frac{1}{2 \gamma} \left( E_0 \frac{E_0 + \sqrt{E_0^2 -
    \epsilon^2}}{\epsilon^2} - \ln \left( \frac{E_0 + \sqrt{E_0^2 
    - \epsilon^2}}{|\epsilon|} \right) - 1 \right)
\end{equation}

If  \eqref{eq: approximation} was exact, $\tau$ would be the time when $E=|\epsilon|$, 
corresponding to the bifurcation for the
undamped system, therefore the solution would be valid up until $\tau$. Since 
this is not the case, $\tau$ instead serves as an approximation for the time
at which \eqref{eq: approximation} stops working. In fact, 
the approximation stops working a small amount of time before $\tau$. This is
expected since \eqref{eq: 2siteshift} does not make sense for 
$E > \sqrt{\epsilon^2 - \tfrac{\gamma^2}{4}}$ due to the domain of the
$\arcsin$.

Figure \ref{fig:2site approximation} shows the effectiveness of the approximation.
$E_0$ was chosen so that it coincided with the initial energy of the graphs in 
the previous sections, leading to $\tau=4700.8$. The energy was graphed for longer
than $\tau$ to identify the behaviour after it leaves the breather, while the
phase difference was graphed for less than $\tau$ since the expression for the
approximation explodes due to the singularity in the $\arcsin$. After this,
the behavior of the
solution is governed either by \eqref{eq: expdecay1} or \eqref{eq: expdecay2}.
Note that by approximating the energy of the solution by the expression for
the energy of the breather in \eqref{eq:two-site_evolution} we are assuming
that this approximation moves along the family of breathers and from Figure 
\ref{fig:2site approximation} we again see that this scenario accurately
reproduces the observed numerical behavior of the solution of the system.

\section{Summary and Conclusions}

In this paper we have derived a new family of periodic solutions of the damped and driven discrete nonlinear Schr\"odinger equation.  We have derived approximations to these solutions and analyzed their stability.  We have also proposed an explanation for the appearance of very long-lived metastable states in the phase space of weakly damped lattice systems in which the trajectory is attracted to a cylinder of such breathers and then drifts very slowly along the cylinder of breathers until this family of solutions disappears.  We suspect that the disappearance occurs near a bifurcation of the undamped system. Finally, we have shown that certain aspects of the metastable states are better approximated by these damped and driven breathers than by the breather solutions of the undamped lattice system. 

\vskip.1in

\noindent
{\bf Acknowledgements:}  The research of CEW was supported in part by
NSF grant DMS-18133.
CEW also acknowledges many helpful conversations about breather 
solutions in the 
dNLS equation with J.-P. Eckmann.
DAC is thankful for the support from Boston University's
Undergraduate Research Opportunities Program.

{

\appendices
\section{Proof of Theorem \ref{th:last_phase}}\label{ap:IFT}

In this appendix, we apply the implicit function theorem to establish the existence of 
the breather solutions discussed in Theorem \ref{th:last_phase}.
We begin by defining a function whose components are the equations for $\dot{p}_j$ and $\dot{q}_j$.  
Thus, we set
\begin{eqnarray}\label{eq:f2-Nm}
&& f_j(p_2,p_3,\dots,p_N,q_2,q_3,\dots,q_N;\epsilon,\omega,\gamma,p_1,q_1) = \\ \nonumber
&& \qquad \qquad \qquad \qquad = \epsilon(q_{j+1}+q_{j-1}-2 q_j) +\omega q_j -(p_j^2+q_j^2)q_j  \\
&& f_{j+N}(p_2,p_3,\dots,p_N,q_2,q_3,\dots,q_N;\epsilon,\omega,\gamma,p_1,q_1) = \\ \nonumber
&& \qquad \qquad  \qquad \qquad = -\epsilon(p_{j+1}+p_{j-1}-2 p_j)  - \omega p_j +(p_j^2+q_j^2)p_j\ ,
\end{eqnarray}
for $j=2,3, \dots N-1$, and 
\begin{eqnarray}\label{eq:fN}
&& f_N(p_2,p_3,\dots,p_N,q_2,q_3,\dots,q_N;\epsilon,\omega,\gamma,p_1,q_1) = \\ \nonumber
&& \qquad \qquad \qquad \qquad =  \epsilon(q_{N-1}-q_N) +\omega q_N -(p_N^2+q_N^2)q_N - \gamma p_N \\
&& f_{2N}(p_2,p_3,\dots,p_N,q_2,q_3,\dots,q_N;\epsilon,\omega,\gamma,p_1,q_1) = \\ \nonumber
&& \qquad \qquad  \qquad \qquad = -\epsilon(p_{N-1}-p_N)  - \omega p_N+(p_N^2+q_N^2)p_N - \gamma q_N .
\end{eqnarray}

Defining
\begin{eqnarray} \nonumber
&& F(p_2,p_3,\dots,p_N,q_2,q_3,\dots,q_N;\epsilon,\omega,\gamma,p_1,q_1) = \\
&& \qquad \qquad =
(f_2, f_3, \dots , f_N, f_{N+2}, f_{N+3}, \dots,  f_{2N})\ ,
\end{eqnarray}
we see that $F: \real^{2(N-1)} \times \real^5 \to \real^{2(N-1)}$ and that
$F(0,0,\dots,0;0,1,0,1,0)=0$.
Computing the Jacobian at this fixed point gives:
\begin{equation}
D_{(\tilde{p},\tilde{q})} F = \left( \begin{array}{cc} 0 & 1 \\ -1 & 0 \end{array} \right)
\end{equation}
where $\tilde{p} = (p_2,\dots , p_N)$, $\tilde{q}=(q_2,\dots , q_N)$, and 
where $0$ represents and $(N-1)\times (N-1)$ matrix of zeros and $1$ represents an $(N-1) \times (N-1)$
dimensional identity matrix.

Thus, by the Implicit Function Theorem, for every $\epsilon, \gamma, q_1$ sufficiently close to zero
and every $\omega$ and $p_1$ sufficiently close to $1$, there exists (a unique)
$p_2^*, \dots , p_N^*, q_2^*, \dots, q_N^*$ such that
\begin{equation}
F(p_2^*,p_3^*,\dots,p_N^*,q_2^*,q_3^*,\dots,q_N^*;\epsilon,\omega,\gamma,p_1,q_1) = 0\ .
\end{equation}
Furthermore, the solution $(p_2^*, \dots , p_N^*, q_2^*, \dots, q_N^*)$ depends
smoothly (in fact, analytically) on $(\epsilon,\omega,\gamma,p_1,q_1)$.

Finally, consider the two equations for $\dot{p}_1$ and $\dot{q}_1$. In order to have a fixed point,  we need
\begin{eqnarray}\label{eq:f1}
&& f_1(\epsilon, \omega, \beta, p_1, q_1, p_2, q_2) = \\ \nonumber
&& \qquad =
\epsilon(q_2-q_1) + \omega q_1 -(q_1^2+p_1^2) q_1 + \beta p_1 = 0 \\
&&f_{N+1}(\epsilon, \omega, \beta, p_1, q_1, p_2, q_2) =
\\ \nonumber
&& \qquad 
= -\epsilon(p_2-p_1) - \omega p_1 +(q_1^2+p_1^2) p_1 + \beta q_1 = 0
\end{eqnarray}

By rotational invariance, we choose $q_1=0$.  Then, inserting the solutions
$q_2^*$ and $p_2^*$ from above, from the requirement that $f_1=0$, we see
that we must have
\begin{equation}\label{eq:beta}
\beta p_1 = - \epsilon q_2^* \ ,
\end{equation}
while the requirement that $f_{N+1} =0$ implies
\begin{equation}\label{eq:p1}
-\epsilon(p_2^*-p_1) - \omega p_1 + p_1^3 = 0\ .
\end{equation}
Using the implicit function theorem, and the fact that $p_2^*$ depends smoothly on 
$p_1$, we see that there exists a solution $p_1^*$ of \eqref{eq:p1}, for $\omega $ near
$1$ and $\gamma$ and $\epsilon$ sufficiently small which satisfies
$(p_1^*)^2= \omega +{\cal{O}}(\epsilon)$.  Inserting this into
\eqref{eq:beta}, we see that there is a unique value of $\beta$ 
(depending smoothly on $\epsilon$, $\gamma$, and $\omega$), such that
we have a damped and driven breather near the $(1,0,0, \dots , 0)$ breather, for each
value of $\omega$ near one and $\epsilon$ and $\gamma$ sufficiently small.

\section{Detailed analysis of the two site system}\label{ap: two-site analysis}

The simplicity of the two-site system allows for a more global investigation of the phase space.   Once again, the  analysis is better suited to
using the energy-phase equations rather than expressing the system in the $p,q$
coordinates, and the energy-phase equations for the two-site system give the simple system of three ODEs:

\begin{eqnarray}\label{eq:2sitedd}
\dot{E}_1 &=& 2 \epsilon \sqrt{E_1 E_2} \sin \psi + 2 \beta E_1 \label{eq:2site}\\ \nonumber
\dot{E}_2 &=& - 2 \epsilon \sqrt{E_1 E_2} \sin \psi - 2 \gamma E_2\\ \nonumber
\dot{\psi} &=& 2 ( E_2 - E_1 ) \left( 1 + \epsilon \frac{\cos \psi}{2 \sqrt{E_1 E_2}} \right)
\end{eqnarray}

For the undriven-undamped, two-site system, (i.e. $\gamma=\beta=0$),  all breather solutions
were found in \cite{Eilbeck:1985} and in our
notation they take the form:
 
\begin{eqnarray}\label{eq:Eilbeck}
    E_1^{(1)} = \tfrac{1}{2} E \quad&     E_2^{(1)} = \tfrac{1}{2} E &\quad    \psi^{(1)} = 0 \label{eq:e1}\\
    E_1^{(2)} = \tfrac{1}{2} E \quad&     E_2^{(2)} = \tfrac{1}{2} E &\quad    \psi^{(2)} = \pi \label{eq:e2}\\
    E_1^{(3)} = \tfrac{1}{2} \left( E + \sqrt{E^2 - \epsilon^2} \right) \quad& E_2^{(3)} = \tfrac{1}{2} \left( E - \sqrt{E^2 - \epsilon^2} \right) &\quad \psi^{(3)} = \pi\\
    E_1^{(4)} = \tfrac{1}{2} \left( E - \sqrt{E^2 - \epsilon^2} \right) \quad& E_2^{(4)} = \tfrac{1}{2} \left( E + \sqrt{E^2 - \epsilon^2} \right) &\quad  \psi^{(4)} = \pi
\end{eqnarray}
 
when $\epsilon>0$ with similar expressions for $\epsilon<0$.  For each of these breathers, we can calculate the angular
frequency $\omega$ from $E$. The values are $\omega^{(1)} = E, \omega^{(2)} = E + 2 \epsilon, \omega^{(3)} = \omega^{(4)} = 2E + \epsilon$. The third breather is the one close to the metastable states discussed so far.

These breathers continue to the damped and driven setting , and again, one can derive explicit formulas for them by calculating the fixed points of \eqref{eq:2sitedd} - for example, 
if $0 < \beta << \gamma$, and $\beta \gamma < \epsilon^2$, the family of breathers $E^{(3)}$, continues to a family of the form described by \eqref{eq:dd2}
(Note that to compare these formulas with those
derived for the general $N$-site damped and driven breather, one should first fix $\gamma$ and $\epsilon$ and then choose $\beta$ 
so that $E^{(3)}_1 = 1/2$, since those breathers were chosen to be close to the breather with $p_1 =1$ and all other components zero.)

If we linearize energy-phase system system \eqref{eq:2site} about this breather, the Jacobian matrix takes the form:

\begin{equation}
    \left( \begin{array}{ccc}
    \beta & - \gamma &- ( \varepsilon^2 - \gamma \beta )  \\
    \beta & - \gamma & ( \varepsilon^2 - \gamma \beta ) \\
    - 1 + \frac{\beta}{\gamma} & 1 - \frac{\gamma}{\beta} & \beta - \gamma
    \end{array} \right)
\end{equation}

Using Mathematica to derive series approximations for the solutions of the characteristic polynomial, 
we find the eigenvalues of the linearization are  the eigenvalue $\lambda_2$, given by \eqref{eq:l22s} and:

\begin{eqnarray}
\lambda_3 &=& -\gamma \left( 1 + \left( \frac{\epsilon}{\omega} \right)^2 +\mathcal{O}(\epsilon^4) \right) + i \left( \omega + \mathcal{O}(\epsilon^2) \right) \label{eq:eig3}\\
\lambda_4 &=& -\gamma \left( 1 + \left( \frac{\epsilon}{\omega} \right)^2 +\mathcal{O}(\epsilon^4) \right) - i \left( \omega + \mathcal{O}(\epsilon^2) \right) \label{eq:eig4}
\end{eqnarray}

These two eigenvalues are close to $\pm i \omega$, but with strictly negative real
part, proportional to $\gamma$,  as expected. An important fact to notice is that if we instead find the linearization about this fixed point using the $p,q$ coordinates, we obtain the same eigenvalues $\lambda_2,\lambda_3,\lambda_4$, with an added $0$ 
eigenvalue corresponding to the invariance with respect to rotations about the cylinder of breathers.

For $|\epsilon|>>\gamma>0$ but $\beta=0$, the breathers disappear but they can be used to help find approximate
solutions of the damped problem, and to understand its metastable behavior.
 For example, if we make the {\it Ansatz} that $E_1(t)=E_2(t)$ for all time, and
substitute this into the equations of motion we find exact solutions: 
 
\begin{eqnarray}
E_1(t) = \tfrac{1}{2} E_0 e^{-\gamma t} \quad & E_2(t) = \tfrac{1}{2} E_0 e^{-\gamma t} & \quad \psi(t) = \ -\arcsin\left( \frac{\gamma}{2\epsilon} \right) \label{eq: expdecay1}\\ 
E_1(t) = \tfrac{1}{2} E_0 e^{-\gamma t} \quad & E_2(t) = \tfrac{1}{2} E_0 e^{-\gamma t} & \quad \psi(t) = -\pi +\arcsin\left( \frac{\gamma}{2\epsilon} \right) \label{eq: expdecay2}
\end{eqnarray}
 
where $E_0=E_1(0)+E_2(0)$ is the initial total energy, and $\epsilon >0$.  One has similar expressions when $\epsilon <0$. These are the states corresponding to \eqref{eq:e1} and \eqref{eq:e2} respectively, with only the first of these being stable.

We can use a similar argument to approximate the effects of initial conditions close
to the solution $E^{(3)}$ of \eqref{eq:Eilbeck}.
  First, differentiate $\tfrac{d}{dt} E_1^{(3)}$ and $\tfrac{d}{dt} E_2^{(3)}$:
\begin{eqnarray} \nonumber
\frac{d}{dt} E_1^{(3)} &=& \frac{1}{2} \dot{E} \left( 1 + \frac{E}{\sqrt{E^2 - \epsilon^2}} \right) \\ \nonumber \frac{d}{dt} E_2^{(3)} &=&  \frac{1}{2} \dot{E} \left( 1 - \frac{E}{\sqrt{E^2 - \epsilon^2}} \right)
\end{eqnarray}
Using $\dot{E} = \dot{E}_1 + \dot{E_2} = - 2 \gamma E_2$, and considering $E_1 = E_1^{(3)}, E_2 = E_2^{(3)}$, we set:
\begin{equation}
    \tfrac{d}{dt} E_1^{(3)} - \tfrac{d}{dt} E_2^{(3)} = \dot{E}_1 - \dot{E}_2=-2\epsilon\sqrt{E_1E_2}\sin\psi-(2\epsilon\sqrt{E_1E_2}-2\gamma E_2)
\end{equation} 

Solving for $\sin \psi$ we get \eqref{eq: 2siteshift}. However, if we calculate $\dot\psi$ from this expression, it will not match the expression in the equations of motion, meaning that the solution is not exact. This is expected due to the metastability of the state, the solution cannot be described by these equations for all time. Numerically, the above expression provides a very good approximation for long time.
Recall that for this breather $\omega^{(3)} = 2E + \epsilon$, so calculating
the leading order in $\psi$ we find $\psi = \tfrac{\gamma}{\omega} + \dots$, which is
identical to the value
of the twist in the damped and driven breather from Theorem \ref{th:twist_approx}. This is further evidence of the fact that
the damped and driven breathers provide a better approximation to
the metastable states than do the undamped breathers.

We also obtain an approximation to the time evolution of the energy through:
\begin{equation}\label{eq:two-site_evolution}
    \dot{E} = - 2 \gamma E_2 \approx - 2 \gamma E_2^{(3)} = - \gamma \left( E - \sqrt{E^2 - \epsilon^2} \right)
\end{equation}
If the system starts with initial total energy $E_0$, an approximate expression we simply solve the integral:

\begin{equation}
    \int_{E_0}^{E(t)}\frac{dx}{x-\sqrt{x^2-\epsilon^2}}=-\gamma t
\end{equation}

If we set the upper limit to be $\epsilon$ and the right hand side to be $-\gamma\tau$, we obtain \eqref{eq:tau}, using the indefinite integral:
\begin{equation}\label{eq:indef}
    \int\frac{dx}{x-\sqrt{x^2-\epsilon^2}}=\frac12\left[x\frac{x+\sqrt{x^2-\epsilon^2}}{\epsilon^2}-\ln\left(x+\sqrt{x^2-\epsilon^2}\right)\right]
\end{equation}
To obtain \eqref{eq: approximation}, it is worth defining:
\begin{equation}
    f(t)=E(t)+\sqrt{E(t)^2-\epsilon^2}
\end{equation}
Then, evaluating the indefinite integral \eqref{eq:indef} at $E(t),E_0$, will lead to an algebraic equation in $f(t)$:
\begin{equation}
    4\gamma t-2\frac{E_0}{\epsilon^2}+1=-\left(\frac{f(t)}{\epsilon}\right)^2+\ln\left(\left(\frac{f(t)}{f(0)}\right)^2\right)
\end{equation}
Which can be shown to solve to:
\begin{equation}
    f(t)=\epsilon\sqrt{-W_{-1}(e^{-1+4\gamma(t-\tau)})}
\end{equation}
From where \eqref{eq: approximation} is derived, and $\tau$ defined as in \eqref{eq:tau}. To get that the branch of the product logarithm, we take the limit:
\begin{equation}
    \lim_{t\to0}E(t)=E_0
\end{equation}

}

\bibliographystyle{unsrt}

\end{document}